\documentclass{ifacconf}

\usepackage{graphicx}      % include this line if your document contains figures
\usepackage{amsmath,amssymb}
\usepackage{xcolor}
\usepackage[inline]{enumitem}
\usepackage{natbib}        % required for bibliography
\usepackage{graphicx,import}
\usepackage{transparent}

\DeclareMathOperator*{\grad}{ \text{\textbf{grad}}}

\renewcommand{\div}{\operatorname{div}}

\newtheorem{remark}{Remark}

\begin{document}
\begin{frontmatter}

\title{Scattering-Passive Structure-Preserving Finite Element Method for the Boundary Controlled Transport Equation with\\ a Moving Mesh
%\thanksref{
} 

\thanks[footnoteinfo]{Corresponding author: Jesus-Pablo Toledo-Zucco. E-mail: jtoledoz@onera.fr}

\author[First,Second]{Jesus-Pablo Toledo-Zucco},
\author[Second]{Denis Matignon},
\author[First,Second]{Charles Poussot-Vassal},

\address[First]{ONERA-DTIS, Université de Toulouse, France}
\address[Second]{Fédération ENAC ISAE-SUPAERO ONERA, Université de Toulouse, France.}
%\address[Second]{AIRBUS Operation SAS, Loads \& Aeroelastics (e-mails: robin.vernay@airbus.com, fabio.vetrano@airbus.com)}
%\address[Third]{ONERA, (e-mail: author@snu.ac.kr)}
%\address[Fourth]{ONERA, (e-mail: author@snu.ac.kr)}

\begin{abstract}                % Abstract of not more than 250 words.
A structure-preserving Finite Element Method (FEM) for the transport equation in one- and two-dimensional domains is presented. This Distributed Parameter System (DPS) has  non-collocated boundary control and observation, and reveals a {\it scattering-energy preserving} structure. We show that the discretized model preserves the aforementioned structure from the original infinite-dimensional system. 
Moreover, we analyse the case of moving meshes for the one-dimensional case. The moving mesh requires less states than the fixed one to produce  solutions with a comparable accuracy, and it can also reduce the overshoot and oscillations of Gibbs phenomenon produced when using the FEM.
% (\textcolor{red}{Gilber efect ? I dont remember the name}). 
Numerical simulations are provided for the case of a one-dimensional transport equation with fixed and moving meshes.
\end{abstract}


\end{frontmatter}
%===============================================================================

\section{Introduction}
Simulation and control design of distributed parameter systems (DPS) described by Partial Differential Equations (PDEs) requires spatial and time discretizations methods. The field of computational methods seeks to develop high-fidelity models with a low computational effort in such a way that they are implementable on available computers.   

One way to reduce computational efforts when multi-physical systems are involved is to use a modular approach. A modular approach allows to represent the overall system as the interconnection of subsystems, which are more simple to deal with. The port-Hamiltonian (pH) approach, introduced in \citep{Maschke1992ConferencePort}, aims at representing complex multi-physical systems as the interconnection between sub-systems by mean of conservation laws. This approach has shown a great success during the last two decades concerning modelling and control of infinite-dimensional systems \citep{Rashad2020JournalTwenty}.

For boundary controlled and observed infinite-dimensional systems with an impedance-passive structure, the pH approach has been widely used for the structure-preserving discretization. In the literature, one can find structure-preserving discretization methods using discrete exterior calculus \citep{Seslija2012JournalDiscrete}, mixed finite-elements \citep{Golo2004JournalHamiltonian}, finite volume \citep{Kotyczka2016ConferenceFinite}, finite-differences \citep{Trenchant2018JournalFinite}, among others. Recently, the Partitioned Finite Element
Method (PFEM) \citep{CardosoRibeiro2018ConferenceStructure}, \citep{Serhani2019ConferencePartitioned}, \citep{Cardoso2020JournalPartitioned} has shown a wide number of applications in one, two and three dimensional domains \citep{Brugnoli2019JournalPort}, \citep{Haine2022ConferenceStructure}. The methods mentioned above are mainly focused in {\it impedance-energy preserving} systems and not much attention has been paid to {\it scattering-energy preserving} systems.

In this paper, the application of the FEM to {\it scattering systems} is presented for the cases of the one- and two-dimensional transport equations. We follow the approach proposed in \citep{Cardoso2020JournalPartitioned} since it has shown to be adaptable to deal with different spatial domains, it uses a single mesh, it provides sparse matrices and it allows distributed parameters. Similar as in the PFEM, in which the {\it impedance-energy preserving} structure is preserved at the discrete level, in the FEM proposed in this paper, the {\it scattering-energy preserving} structure is preserved at the discrete level.

A second contribution of this paper is the extension of the structure-preserving FEM to deal with moving meshes. It is known that the efficiency of numerical methods for solving PDEs can be improved concentrating the nodes in the areas of rapid variations of the solution \citep{Huang1994JournalMoving}, \citep{Huang2002JournalAdaptive}, \citep{Aydougdu2019JournalData}. In this paper, some first steps towards the structure-preserving  FEM using moving meshes is presented for the case of the one-dimensional transport equation. The discretized model is a time-variant finite-dimensional system, in which the time-variying matrices depend on the mesh dynamic through sparse matrices as well. We show that the {\it scattering} structure is also preserved on the finite-dimensional model when the velocity density is constant in space. 

In Section \ref{sec:TransportEquation}, we recall the transport equation in one- and two-dimensional domains, showing it {\it scattering} structure with respect to a quadratic Hamiltonian. In Section \ref{sec:PFEM}, we develop the structure-preserving discretization scheme for both cases, showing that the discretized models are also {\it scattering-energy systems}. In Section \ref{sec:MovingMesh}, we present some first steps towards moving meshes using the toy example of the one-dimensional transport equation. In Section \ref{sec:Simulations}, we show some simulations on the one-dimensional transport equation with fix and moving meshes. Finally, in Section \ref{sec:conclusions} we give some conclusions and perspectives.

\section{Transport Equation}\label{sec:TransportEquation}
\subsection{One-dimensional transport equation}
We consider the following transport equation in one-dimensional space:
\begin{equation}\label{TransportEquation}
\begin{cases}
\dot{x}(\zeta,t) +c \dfrac{\partial e}{\partial \zeta} (\zeta,t)=0, \\ x(\zeta,0) = x_0(\zeta), \\
e(\zeta,t) = \mathcal{H}(\zeta)x(\zeta,t),
\end{cases}
\end{equation}
in which $ \zeta \in [a,b]=:\Omega$ represents the space, $ t \geq 0$ the time, $x(\zeta,t) \in \mathbb{R}$ the state, $e(\cdot,t) \in H^1(\Omega)$ the co-energy variable, $c \in \lbrace -1, +1 \rbrace$ the velocity direction, $\mathcal{H}(\zeta) >0 $ the Hamiltonian density, and $x_0(\zeta)$ the initial condition. The  input $u(t)$ and output $y(t)$ are chosen at the spatial boundaries and they depend on the velocity direction as follows:
\begin{equation}\label{InputOutput}
c = -1: \begin{cases}
e(b,t) = u(t), \\
e(a,t) = y(t),
\end{cases} \,c = +1 :\begin{cases}
e(a,t) = u(t), \\
e(b,t) = y(t),
\end{cases}
\end{equation}
The Hamiltonian associated to \eqref{TransportEquation} is defined as
\begin{equation}\label{Hamiltonian}
H(t) := \dfrac{1}{2}\int_a^b \mathcal{H}(\zeta) x(\zeta,t) ^2 d\zeta=\dfrac{1}{2}\int_a^b  x(\zeta,t)  e(\zeta,t) d\zeta,
\end{equation}
and one can show that the following balance equation is satisfied:
\begin{equation}\label{Balance}
\dot{H}(t) = \dfrac{1}{2} \left(  u(t)^2 - y(t) ^2 \right).
\end{equation}
This balance guarantees that the system \eqref{TransportEquation}-\eqref{InputOutput} is {\it scattering energy preserving}.

Note that, when $\mathcal{H}(\zeta)=\mathcal{H}_0$, then the analytic solution is known as $y(t) = u(t-\tau)$, with $\tau=\frac{\ell}{\mathcal{H}_0}$ and $\ell = b-a$. This is a simple delay system.

\subsection{Two-dimensional transport equation}\label{subsec:2D}

Now, a two-dimensional transport equation is considered. The PDE is:
\begin{equation}\label{Eq:TE2D}
    \dfrac{\partial x}{\partial t}(\zeta_1,\zeta_2,t) + \div(\bold{c}(\zeta_1,\zeta_2)x(\zeta_1,\zeta_2,t)) = 0 
\end{equation}
in which $(\zeta_1,\zeta_2) \in \Omega$ is the space with $\Omega$ a bounded domain, $t\geq 0$ is the time, $x$ is the state variable, and $\bold{c}(\zeta_1,\zeta_2) = \begin{bmatrix}
    c_1(\zeta_1,\zeta_2) \\c_2(\zeta_1,\zeta_2)
\end{bmatrix}$ is the velocity distribution vector. We define the boundary of $\Omega$ as $\Gamma$ and the normal vector $\bold{n}(\zeta_1,\zeta_2)= \begin{bmatrix}
    n_1(\zeta_1,\zeta_2) \\n_2(\zeta_1,\zeta_2)
\end{bmatrix}$ which is perpendicular to $\Gamma$ (pointing outside of $\Omega$). We define $\Gamma_i$ and $\Gamma_u$ as the sets of the input and output, respectively. These subsets are defined as:
\begin{align}
    \Gamma_i &= \lbrace (\zeta_1,\zeta_2) \in \Gamma, \quad \bold{c}(\zeta_1,\zeta_2)^\top \bold{n}(\zeta_1,\zeta_2)<0 \rbrace , \\
    \Gamma_o &= \lbrace (\zeta_1,\zeta_2) \in \Gamma, \quad \bold{c}(\zeta_1,\zeta_2)^\top \bold{n}(\zeta_1,\zeta_2)>0 \rbrace .
\end{align}
These subsets are such that $\Gamma = \Gamma_i \cup \Gamma_o$, and the input and output are:
\begin{equation}\label{InputOutput2D}
\begin{split}
    u (\zeta_1,\zeta_2,t) &= x(\zeta_1,\zeta_2,t), \quad (\zeta_1,\zeta_2) \in \Gamma_i,\\
    y (\zeta_1,\zeta_2,t) &= x(\zeta_1,\zeta_2,t), \quad (\zeta_1,\zeta_2) \in \Gamma_o.
\end{split}
\end{equation}
The Hamiltonian associated to \eqref{Eq:TE2D} is defined as:
\begin{equation}\label{Eq:Hamiltonian2D}
    H(t) : = \dfrac{1}{2}\iint _\Omega x(\zeta_1,\zeta_2,t)^2 d\zeta,
\end{equation}
with $\zeta = (\zeta_1,\zeta_2)$. In the following, for simplicity, we omit the time and spatial dependency. The time derivative of the Hamiltonian is computed as follows:
\begin{equation}\label{Eq:Hdot_2D}
\begin{split}
    \dot{H} &= \iint _\Omega x \dfrac{\partial x}{\partial t} d\zeta = -\iint _\Omega x \, \div (\bold{c}x) d\zeta, \\
    &= \iint _\Omega x \, \bold{grad (x)}^\top \bold{c} \,  d\zeta - \int_\Gamma x^2 \bold{c}^\top \bold{n} d\zeta,
\end{split}
\end{equation}
in which we have applied integration by part using Stokes theorem.

Now, we use the divergence rule for a product:
\begin{equation}\label{Eq:divrule}
    \div(\bold{c}\,x) = \bold{c} ^\top \, \bold{grad(x)} + x\, \div(\bold{c})
\end{equation}
to highlight the {\it scattering energy preserving} structure of the two-dimensional transport equation. To this end, we multiply \eqref{Eq:divrule} by $x$ as follows:
\begin{equation}\label{Eq:divrule2}
    x\, \div(\bold{c}\,x) = x\, \bold{c}^\top \, \bold{grad(x)} + x^2\, \div(\bold{c})
\end{equation}
implying the following:
\begin{equation}
\begin{split}
    \dot{H} &= -\iint _\Omega x \, \div (\bold{c}x) d\zeta \\
    &=  -\iint _\Omega \left[ x\, \bold{c}^\top \, \bold{grad(x)} + x^2\, \div(\bold{c}) \right]d\zeta, \\
    &=  -\iint _\Omega \left[ x\,  \, \bold{grad(x)}^\top \bold{c} + x^2\, \div(\bold{c}) \right] d\zeta. 
\end{split}
\end{equation}
From the previous equation, we use the following identity:
\begin{equation}
    \iint _\Omega  x\,  \, \bold{grad(x)}^\top \bold{c} \,d\zeta = -\dot{H}  - \iint _\Omega   x^2\, \div(\bold{c}) d\zeta
\end{equation}
and replace it in \eqref{Eq:Hdot_2D}
\begin{equation}
\begin{split}
    \dot{H} &= \iint _\Omega x\, \bold{grad (x)}^\top c \,  d\zeta - \int_\Gamma x^2 \bold{c}^\top \bold{n} d\zeta, \\
    &= -\dot{H}  - \iint _\Omega   x^2\, \div(\bold{c}) d\zeta - \int_\Gamma x^2 \bold{c}^\top \bold{n} d\zeta,
\end{split}
\end{equation}
implying
\begin{equation}\label{Eq:Balance2D}
    \dot{H} = -\dfrac{1}{2} \iint _\Omega   x^2\, \div(\bold{c}) d\zeta - \dfrac{1}{2}\int_\Gamma x^2 \bold{c}^\top \bold{n} d\zeta.
\end{equation}
Notice that, if $\div(\bold{c}) = 0$, the system is {\it scattering energy preserving} with the following balance:
\begin{equation}
\begin{split}
    \dot{H} 
   &=  \dfrac{1}{2}\int_{\Gamma_i} u ^2 \,|\bold{c}^\top \bold{n_i}| d\zeta - \dfrac{1}{2} \int_{\Gamma_o} y ^2 \,\bold{c}^\top \bold{n_o} d\zeta, \\
\end{split}
\end{equation}
with $\bold{n_i}$ and $\bold{n_o}$ the normal vectors associated to $\Gamma_i$ and $\Gamma_o$, respectively, and $|\bold{c}^\top \bold{n_i}|$ and $\bold{c}^\top \bold{n_o} $ a density on the surface $\Gamma$.

\subsubsection{Example} 
\begin{figure}[!h]
\begin{center}
\includegraphics[width=0.28\textwidth]{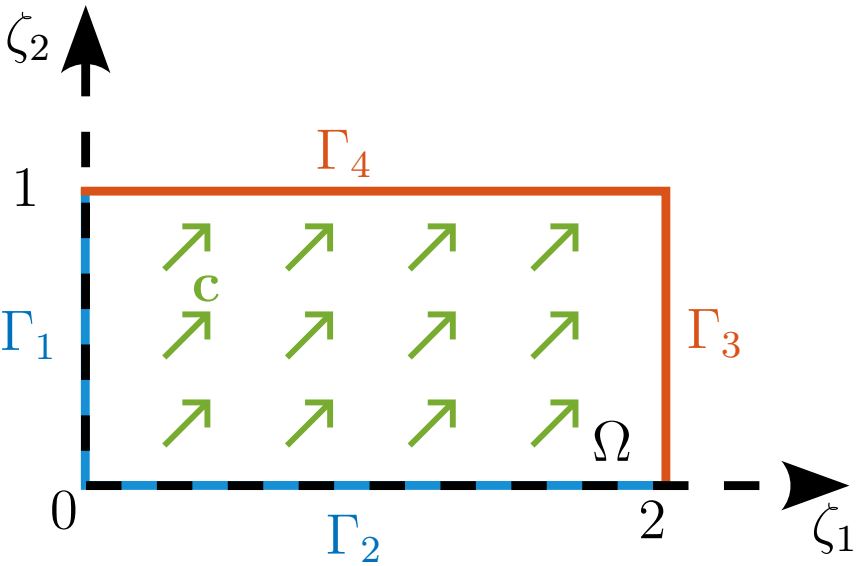}  
\caption{Spatial domain and boundary partitions.} 
\label{Fig:SpatialDomain2D}                      \end{center}
\end{figure}

Consider the case of a rectangular domain $\Omega = [0,2] \times [0,1]$ and $\bold{c} = \begin{bmatrix}
    1 \\ 1
\end{bmatrix}$ for all $(\zeta_1,\zeta_2) \in \Omega$. We define the following boundary domains:
\begin{itemize}
    \item $\Gamma_1 = \lbrace \zeta_1 = 0, \quad \quad  \zeta_2 \in [0,1] \rbrace$
    \item $\Gamma_2 = \lbrace \zeta_1 \in [0,2], \;\;\zeta_2 =0 \rbrace$
    \item $\Gamma_3 = \lbrace \zeta_1 = 2, \quad \quad  \zeta_2 \in [0,1] \rbrace$
    \item $\Gamma_4 = \lbrace \zeta_1 \in [0,2],\;\; \zeta_2 =1 \rbrace$
\end{itemize}
for which the normal vectors are respectively:
\begin{equation}
    \bold{n_1} = \begin{bmatrix}
        -1 \\0
    \end{bmatrix}, \; \bold{n_2} = \begin{bmatrix}
        0 \\-1
    \end{bmatrix},\; \bold{n_3} = \begin{bmatrix}
        1 \\0
    \end{bmatrix},\; \bold{n_4} = \begin{bmatrix}
        0 \\1
    \end{bmatrix}. 
\end{equation}
and since $\bold{c}^\top \bold{n_1} = -1$, $\bold{c}^\top \bold{n_2} = -1$, $\bold{c}^\top \bold{n_3} = 1$, $\bold{c}^\top \bold{n_4} = 1$, 
$\Gamma_i = \Gamma_1 \cup \Gamma_2$ and $\Gamma_o = \Gamma_3 \cup \Gamma_4$. In Fig. \ref{Fig:SpatialDomain2D}, we show an schematic of the spatial domain and the boundary partitions.

\vspace{5mm}
%\subsection*{Contributions of this paper}
In the following section, we show that the proposed discretization scheme preserves the energy properties \eqref{Hamiltonian}-\eqref{Balance} for the one-dimensional case and \eqref{Eq:Hamiltonian2D}-\eqref{Eq:Balance2D} for the two-dimensional case. 

%\subsection{Weak Formulation}
%Let us consider the following test function $v(\zeta,t) \in H^1([a,b],\mathbb{R}_+)$. We write the following weak form by multiplying the PDE in \eqref{TransportEquation} by $v(\zeta,t)$ and integrating over the spatial domain:
%\begin{equation}\label{WeakForm}
%\begin{split}
%\int_a^b v(\zeta,t) \dot{x}(\zeta,t) d\zeta &= \int_a^b v(\zeta,t) \dfrac{\partial x}{\partial \zeta} (\zeta,t) d\zeta, \\
%\end{split}
%\end{equation}
%The main difference with respect to the traditional PFEM approach is that we consider time-dependent test functions. 

\section{Finite Element Method}\label{sec:PFEM}

In this section, the FEM is applied to the one- and two dimensional transport equations introduced previously. We show that in both cases, the discretized Hamiltonian mimics the Hamiltonian of the infinite-dimensional system, implying that the discretized models are also {\it scattering energy systems}.

\subsection{One-dimensional FEM for the transport equation}\label{subsec:1Dcase}
We approximate the state and co-energy variable using the following projections:
\begin{equation}\label{Approximation_fix}
\begin{split}
x(\zeta,t) \approx x^{ap}(\zeta,t) &:= \Phi(\zeta) ^\top x_d(t), \\
e(\zeta,t) \approx e^{ap}(\zeta,t) &:= \Phi(\zeta) ^\top e_d(t).
\end{split}
\end{equation}
The approximated functions $x^{ap}(\zeta,t)$ and $e^{ap}(\zeta,t) \in \mathbb{R}$ for all $\zeta \in \Omega$ and $t \geq 0$ are both projected on the same basis functions 
\begin{equation}
\Phi(\zeta) = \begin{pmatrix}
\phi_1 (\zeta) \\ \vdots \\ \phi_N (\zeta)
\end{pmatrix}
\end{equation}
with $\phi_i(\cdot) \in H^1(\Omega)$ for all $i \in \lbrace 1,\cdots , N \rbrace$. The discretized state and co-energy variables are $x_d(t)$ and $e_d(t) \in \mathbb{R}^{N\times 1}$, respectively.

We apply the following three operations to the equations in \eqref{TransportEquation}: $(i)$ we multiply by $\Phi(\zeta)$ from the right to the left, $(ii)$ integrate over $\zeta \in \Omega$, and $(iii)$ replace the approximation functions from \eqref{Approximation_fix}. For simplicity, from now on, we consider $c = -1$, meaning that the signals travel from the right to the left. The weak form of \eqref{TransportEquation}-\eqref{InputOutput} writes:
\begin{equation}\label{TransportEquationWeak_fix}
\begin{cases}
\int_a^b \Phi(\zeta)\left[\Phi(\zeta)^\top\dot{x}_d(t) -  \dfrac{\partial \Phi}{\partial \zeta} (\zeta)^\top e_d(t)\right] d\zeta =0,\\
\int_a^b \Phi(\zeta) \Phi(\zeta)^\top x_{d}(0) d\zeta = \int_a^b \Phi(\zeta) x_0(\zeta) d\zeta , \\
\int_a^b \Phi(\zeta) \Phi(\zeta)^\top e_d(t) d\zeta = \int_a^b \Phi(\zeta) \mathcal{H}(\zeta) \Phi(\zeta)^\top x_d(t) d\zeta ,\\
y(t) = \Phi(a)^\top e_d(t).
\end{cases}
\end{equation}
which can be equivalently written as:
\begin{equation}\label{Realization_fix}
\begin{cases}
E \dot{x}_d(t) = F e_d(t) + B u(t),\\
 x_d(0) = x_{d0},\\
E e_d(t) = Q x_d(t), \\
y(t) = C e_d(t),
\end{cases}
\end{equation}
with 
\begin{equation}
\begin{split}
E &= \int_a ^b \Phi(\zeta) \Phi(\zeta)^\top d\zeta, \\
F &= -\Phi(a)\Phi(a)^\top -\int_a ^b \dfrac{\partial \Phi}{\partial \zeta}(\zeta) \Phi(\zeta)^\top d\zeta, \\
Q &= \int_a ^b \Phi(\zeta)\mathcal{H}(\zeta) \Phi(\zeta)^\top d\zeta, \\
B & = \Phi(b), \\
C & = \Phi(a), \\
x_{d0} & = E^{-1} \int_a^b \Phi(\zeta)x_0(\zeta)d\zeta.
\end{split}
\end{equation}

%Since $\Phi(\zeta)$ is constant in time, the realization \eqref{Realization_fix} has constant matrices $E$, $F$, $Q$, $B$, and $C$. \textcolor{blue}{Motivate the moving mesh and give some references}.

\subsection*{Discretized Hamiltonian: one-dimensional case}

We replace the approximations \eqref{Approximation_fix} in the Hamiltonian \eqref{Hamiltonian} to obtain the discretized Hamiltonian as follows:
\begin{equation}\label{DiscretizedHamiltonian_fixed}
\begin{split}
H(t) \approx H_d(t) &:= \dfrac{1}{2} \int_a^b x^{ap}(\zeta,t)e^{ap}(\zeta,t) d\zeta, \\
%&= \dfrac{1}{2}\int_a^b \Phi (\zeta)^\top x_d(t) \Phi (\zeta)^\top e_d (t) d\zeta, \\
%&= \dfrac{1}{2}\int_a^b x_d(t)^\top  \Phi (\zeta)\Phi (\zeta)^\top e_d (t) d\zeta, \\
&= \dfrac{1}{2}x_d(t)^\top \left( \int_a^b   \Phi (\zeta)\Phi (\zeta)^\top d\zeta \right) \, e_d (t) , \\
&= \dfrac{1}{2} x_d(t)^\top E e_d (t) =\dfrac{1}{2} x_d(t)^\top Q x_d (t)  .
\end{split}
\end{equation}
The main properties that we use to develop the time derivative of \eqref{DiscretizedHamiltonian_fixed} are:
\begin{align*}
E & = E^\top ,\quad Q = Q^\top,\\
y(t) & = \Phi(a)^\top e_d(t), \\
u(t) & = \Phi(b)^\top e_d(t),\\
E\dot{x}_d(t) &= \int _a^b \Phi(\zeta) \dfrac{\partial \Phi}{\partial \zeta}(\zeta)^\top d\zeta e_d(t).
\end{align*}
For simplicity, we omit time and space dependency. The Hamiltonian balance is:  
\begin{equation}
\begin{split}
\dot{H}_d &=  {x}_d^\top Q \dot{x}_d,\\
&=  {x}_d^\top Q E^{-1} E \dot{x}_d,\\
&=  (Q{x}_d)^\top  E^{-1} (Fe_d + Bu),\\
&=  (E{e}_d)^\top  E^{-1} ((-\Phi_a \Phi_a^\top - \int_a^b \tfrac{\partial \Phi}{\partial \zeta}\Phi ^\top d\zeta)e_d + \Phi_b u),\\
&=  -{e}_d^\top \Phi_a \Phi_a ^\top e_d   -e_d^\top \int_a^b \tfrac{\partial \Phi}{\partial \zeta}\Phi^\top  d\zeta e_d + e_d^\top \Phi_b u,\\
&=  -y^\top y   -e_d^\top \int_a^b \tfrac{\partial \Phi}{\partial \zeta}\Phi^\top d\zeta e_d + u^\top u,\\
&=  -y^2   -e_d^\top \int_a^b \Phi \tfrac{\partial \Phi}{\partial \zeta}^\top d \zeta e_d + u^2,\\
&=  -y^2   -e_d^\top E \dot{x}_d + u^2,\\
&=  -y^2   -(Ee_d)^\top  \dot{x}_d + u^2,\\
&=  -y^2   -(Qx_d)^\top  \dot{x}_d + u^2,\\
&=  -y^2   -x_d^\top Q  \dot{x}_d + u^2,\\
&=  -y^2   -\dot{H}_d + u^2,\\
\Leftrightarrow  \dot{H}_d & = \dfrac{1}{2}(u^2-y^2),
\end{split}
\end{equation}
implying that the discretized system \eqref{Realization_fix} preserves the {\it scattering energy preserving} structure of the infinite-dimensional system.

\subsection{Two-dimensional FEM for the transport equation}
We approximate the state variable using the following projection:
\begin{equation}\label{Approximation_2D}
x(\zeta_1,\zeta_2,t) \approx x^{ap}(\zeta_1,\zeta_2,t) := \Phi(\zeta_1,\zeta_2) ^\top x_d(t), 
\end{equation}
The approximated function $x^{ap}(\zeta_1,\zeta_2,t) \in \mathbb{R}$ for all $(\zeta_1,\zeta_2) \in \Omega$ and $t \geq 0$ is projected on the following basis functions 
\begin{equation}
\Phi(\zeta_1,\zeta_2) = \begin{pmatrix}
\phi_1 (\zeta_1,\zeta_2) \\ \vdots \\ \phi_N (\zeta_1,\zeta_2)
\end{pmatrix}
\end{equation}
with $\phi_i(\cdot,\cdot) \in H^1(\Omega)$ for all $i \in \lbrace 1,\cdots , N \rbrace$. The discretized state variable is $x_d(t) \in \mathbb{R}^{N\times 1}$.

We apply the following three operations to the equations in \eqref{Eq:TE2D}: $(i)$ we multiply by $\Phi(\zeta_1,\zeta_2)$ from the right to the left, $(ii)$ integrate over $\zeta \in \Omega$, and $(iii)$ replace the approximation functions from \eqref{Approximation_2D}. The weak form of \eqref{Eq:TE2D} writes:

\begin{equation}\label{TransportEquationWeak_2D}
  \iint_\Omega \Phi \left[ \Phi^\top \dot{x}_d + \div(\bold{c}\,\Phi ^\top x_d)
 \right]  d\zeta =0
\end{equation}
which can be equivalently written as:
\begin{equation}
    M \dot{x}_d (t) = F_1 x_d (t)
\end{equation}
with 
\begin{align*}
    M &= \iint _ \Omega \Phi \Phi^\top d\zeta , \\
    F_1 &= -\iint _ \Omega \Phi \begin{bmatrix}
        \div(\bold{c}\phi_1) \cdots \div(\bold{c} \phi_N)
    \end{bmatrix} d\zeta 
\end{align*}
or as:
\begin{equation}
    M \dot{x}_d (t) = F_2 x_d (t) - \int_{\Gamma_i}\Phi\,u\,\bold{c}^\top \bold{n_i}\,ds
\end{equation}
with 
\begin{align*}
    F_2 &= \iint _ \Omega \begin{bmatrix}
        \grad(\phi_1)^\top \bold{c}\\
        \vdots \\
        \grad(\phi_N)^\top \bold{c}
    \end{bmatrix} \Phi ^\top  d\zeta - \int _{\Gamma_o} \left( \Phi \Phi^\top \right) \,\bold{c}^\top \bold{n_o}\,ds, 
\end{align*}
in which we have applied integration by parts making appear the input term and keeping the output term inside the internal dynamic through $F_2$.

\subsection*{Discretized Hamiltonian: two-dimensional case}

We replace the approximation \eqref{Approximation_2D} in the Hamiltonian \eqref{Eq:Hamiltonian2D} to obtain the discretized Hamiltonian as follows:
\begin{equation}\label{DiscretizedHamiltonian_2D}
\begin{split}
H(t) \approx H_d(t) &:= \dfrac{1}{2} \iint _\Omega x^{ap}(\zeta,t)x^{ap}(\zeta,t), 
\\&= \dfrac{1}{2} x_d(t)^\top \left( \iint _\Omega \Phi(\zeta)\Phi(\zeta)^\top d\zeta \right)  x_d(t), \\
\end{split}
\end{equation}

Notice that $H_d$ in \eqref{DiscretizedHamiltonian_2D} is defined as $H$ from \eqref{Eq:Hamiltonian2D} but evaluated at $x = x^{ap}$, with $x^{ap}$ from \eqref{Approximation_2D}. Then, following subsection \ref{subsec:2D}, one can also apply the identities \eqref{Eq:divrule}-\eqref{Eq:divrule2} to obtain the following discretized Hamiltonian balance:
\begin{equation}\label{Eq:DiscretizedBalance2D}
    \dot{H} = -\dfrac{1}{2} \iint _\Omega   x^{ap} x^{ap}\, \div(\bold{c})\,d\zeta - \dfrac{1}{2}\int_\Gamma x^{ap} x^{ap} \bold{c}^\top \bold{n}\, ds.
\end{equation}
Then, replacing $x^{ap} = \Phi^\top x_d$ from \eqref{Approximation_2D}, we obtain the following balance:
\begin{equation}\label{DiscretizedHamiltonian_2D}
\begin{split}
{\dot H}_d(t) &= - \dfrac{1}{2} x_d(t)^\top \left( \iint _\Omega \Phi(\zeta) \div(\bold{c})  \Phi(\zeta)^\top d\zeta \right)  x_d(t) 
\\&+ \dfrac{1}{2}\int_{\Gamma_i} {u}^2 \,|\bold{c}^\top \bold{n_i}|\,ds - \dfrac{1}{2} \int_{\Gamma_o} {y}^2 \,\bold{c}^\top \bold{n_o}\, ds,
\\
\end{split}
\end{equation}
where we have replaced the inputs and outputs from \eqref{InputOutput2D}. The previous balance equation shows that at the discretized model preserves the {\it scattering-energy preserving} structure.
\begin{remark}
    For a {\em divergence-free} velocity field, the {\it scattering energy preserving} structure is exactly preserved. However, in a non {\em divergence-free} case, a dissipative or accumulative effect is appreciated depending on the sign of $\div(\bold{c})$ (See \eqref{Eq:Balance2D}) The dissipative or accumulative effect is also preserved at the discrete level through the matrix $ Q_c = \iint _\Omega \Phi \div(\bold{c}) \Phi ^\top d\zeta$
\end{remark}

\begin{remark}
    For the numerical implementation of the discretized model of the 2D transport equation, the input and output have to be projected on two different finite-dimensional bases as well. These projections lead to a new energy balance which is now approximated. The same holds for the {\it scattering energy preserving} structure.
\end{remark}

%Note that since $u = \Phi_b x_d$, we have thee following identity
%\begin{equation}
%\begin{split}
%E \dot{x}_d &= Ax_d + Bu, \\
%&= \int_a^b \left[ \Phi  \tfrac{\partial \Phi}{\partial \zeta}^\top - \Phi \dot{\Phi}^\top \right] x_d.  
%\end{split} 
%\end{equation}
%We replace it in the balance equation of $H_d$ 
%\begin{equation}
%\begin{split}
%\dot{H}_d &=- y^2 +u^2 +\tfrac{1}{2}  x_d^\top  \dot{E} {x}_d - 2{x}_d^\top \int _a^b \Phi \dot{\Phi} ^\top d\zeta  x_d   \\
%& \quad - {x}_d^\top E  \dot {x}_d   \\
%&=- y^2 +u^2 + x_d^\top  \int _a^b\Phi \dot{\Phi} ^\top d\zeta {x}_d - 2{x}_d^\top \int _a^b \Phi \dot{\Phi} ^\top d\zeta  x_d  \\ 
%& \quad - {x}_d^\top E  \dot {x}_d , \quad\quad \textcolor{blue}{ \tfrac{1}{2}  x_d^\top  \dot{E} {x}_d= x_d^\top \int_a^b \Phi \dot{\Phi} ^\top d\zeta {x}_d }  \\
%&=- y^2 +u^2  - {x}_d^\top \int _a^b \Phi \dot{\Phi} ^\top d\zeta  x_d   - {x}_d^\top E  \dot {x}_d ,  \\
%&=- y^2 +u^2  - \left( {x}_d^\top \int _a^b \Phi \dot{\Phi} ^\top d\zeta  x_d   + {x}_d^\top E  \dot {x}_d \right),  \\
%&=- y^2 +u^2  - \left( \tfrac{1}{2} {x}_d^\top \dot{E} x_d   + {x}_d^\top E  \dot {x}_d \right),  \\
%&=- y^2 +u^2  - \dot{H}_d ,  \\
%\end{split}
%\end{equation}
%implying
%\begin{equation}
%\dot{H}_d = \dfrac{u^2}{2} - \dfrac{y^2}{2}
%\end{equation}
%which mimics the initial balance in \eqref{Balance}.
%
%\begin{remark}
%The structure is preserved for any time parametrization of the test functions.
%\end{remark}

\section{Moving Mesh for the 1D case}\label{sec:MovingMesh}

Now, we consider a moving mesh. The state and co-energy variables are approximated as follows:
\begin{equation}\label{Approximation}
\begin{split}
x(\zeta,t) \approx x^{ap}(\zeta,t) &:= \Phi(\zeta,t) ^\top x_d(t), \\
e(\zeta,t) \approx e^{ap}(\zeta,t) &:= \Phi(\zeta,t) ^\top e_d(t)
\end{split}
\end{equation}
with $x^{ap}(\zeta,t)$ and $e^{ap}(\zeta,t) \in \mathbb{R}$ for all $\zeta \in \Omega$ and $t \geq 0$,\begin{equation}
\Phi(\zeta,t) = \begin{pmatrix}
\phi_1 (\zeta,t) \\ \vdots \\ \phi_N (\zeta,t)
\end{pmatrix}
\end{equation}
the vector with the approximation basis functions $\phi_i(\cdot,t) \in H^1(\Omega)$ for all $t \geq 0$ and $i \in \lbrace 1,\cdots , N \rbrace$.

The main difference with respect to the previous case is that, in a moving mesh, the time derivative of the approximation function $x^{ap}(\zeta,t)$ has the following two terms:
\begin{equation}\label{Eq:EdotMoving}
\dot{x}^{ap}(\zeta,t) = \dfrac{\partial \Phi}{\partial t}(\zeta,t) ^\top x_d(t) + \Phi(\zeta,t)^\top \dot{x}_d(t) .
\end{equation} 

We apply the following three operations to the transport equation in \eqref{TransportEquation}: $(i)$ we multiply by $\Phi(\zeta,t)$ from the right to the left, $(ii)$ integrate over $\zeta \in \Omega$, and $(iii)$ replace the approximation function from \eqref{Approximation}-\eqref {Eq:EdotMoving}. Considering $c = 1$, the weak form of \eqref{TransportEquation}-\eqref{InputOutput} writes:

\begin{equation}\label{TransportEquation_weak_moving}
\begin{cases}
\int_a^b \Phi(\zeta,t) f(\zeta,t) d\zeta = 0, \\
f = \dfrac{\partial \Phi}{\partial t}(\zeta,t) ^\top x_d(t) + \Phi(\zeta,t)^\top\dot{x}_d(t) -  \dfrac{\partial \Phi}{\partial \zeta} (\zeta,t)^\top e_d(t),\\
\int_a^b \Phi(\zeta,0) \Phi(\zeta,0)^\top x_{d}(0) d\zeta = \int_a^b \Phi(\zeta,0) x_0(\zeta) d\zeta , \\
\int_a^b \Phi(\zeta,t) \Phi(\zeta,t)^\top e_d(t) d\zeta \\
\quad \quad \quad\quad \quad = \int_a^b \Phi(\zeta,t) \mathcal{H}(\zeta) \Phi(\zeta,t)^\top x_d(t) d\zeta ,\\
y(t) = \Phi(a,t)^\top e_d(t).
\end{cases}
\end{equation}
The realization \eqref{TransportEquation_weak_moving} can be equivalently written as:
\begin{equation}
\begin{cases}
E(t) \dot{x}_d(t) = F(t) e_d(t) + B(t) u(t),\\
x_d(0) = x_{d0}, \\
E(t)e_d(t) = Q(t)x_d(t),\\
y(t) = C(t) x_d(t),
\end{cases}
\end{equation}
with
\begin{equation}\label{Eq:MatricesMovingMesh}
\begin{split}
E(t) & = \int_a^b \Phi(\zeta,t) \Phi(\zeta,t)^\top d\zeta, \\
F(t) & = -\Phi(a,t)\Phi(a,t)^\top  -\int _a^b \dfrac{\partial \Phi}{\partial \zeta}(\zeta,t)\Phi(\zeta,t)^\top d\zeta \\
& \quad - \left(\int _a^b  \Phi(\zeta,t) \dfrac{\partial \Phi}{\partial t} (\zeta,t)^\top d\zeta \right)Q(t)^{-1}E(t), \\
Q(t) &= \int_a ^b \Phi(\zeta,t)\mathcal{H}(\zeta) \Phi(\zeta,t)^\top d\zeta, \\
B(t) &= \Phi (b,t), \\
C(t) & = \Phi(a,t)^\top, \\
x_{d0} & = E(0)^{-1} \int_a^b \Phi(\zeta,0)x_0(\zeta)d\zeta.
\end{split}
\end{equation}

\begin{remark}
    Notice that for $\mathcal{H}(\zeta) = \mathcal{H}_0$ (constant in space), the term $Q(t)^{-1}E(t)$ in the expression of $F(t)$ from \eqref{Eq:MatricesMovingMesh} becomes the identity matrix times $\mathcal{H}_0$.
\end{remark}
%Replacing \eqref{TestFunction} in \eqref{WeakForm} and using the fact that \eqref{WeakForm} should be valid for any $v_d$ we obtain:
%\begin{equation}\label{WeakForm2}
%\begin{split}
%\int_a^b \Phi(\zeta,t) \dot{x}(\zeta,t) d\zeta 
%&= [ \Phi(\zeta,t)  x(\zeta,t)]_a^b  \\ &\quad - \int_a^b  \dfrac{\partial \Phi}{\partial \zeta} (\zeta,t) x(\zeta,t)d\zeta.
%\end{split}
%\end{equation}
%
%Replacing the approximation \eqref{Approximation} into \eqref{WeakForm2}, we obtain the following left-hand side term
%\begin{equation}\label{LeftHand}
%\int_a^b \Phi(\zeta,t) \Phi(\zeta,t)^\top d\zeta \dot{x}_d(t) + \int_a^b \Phi(\zeta,t) \dot{\Phi}(\zeta,t)^\top d\zeta {x}_d(t),
%\end{equation}
%and the following right-hand side term:
%\begin{equation}\label{RightHand}
%\left( \left[\Phi(\zeta,t) \Phi(\zeta,t) ^\top \right]_a^b - \int _a^b \dfrac{\partial \Phi}{\partial \zeta}(\zeta,t) \Phi(\zeta,t) d\zeta \right) x_d(t)
%\end{equation}
%Combining both \eqref{RightHand}-\eqref{LeftHand}, one can write the following time-variant realization:

\subsection*{Discretized Hamiltonian using a moving mesh}

In this case, we replace the approximations \eqref{Approximation} in the Hamiltonian \eqref{Hamiltonian} to obtain the discretized Hamiltonian as follows:
\begin{equation}\label{DiscretizedHamiltonian}
\begin{split}
 H_d(t) &:= \dfrac{1}{2} \int_a^b x^{ap}(\zeta,t)e^{ap}(\zeta,t), \\
&= \dfrac{1}{2}\int_a^b \Phi (\zeta,t)^\top x_d(t) \Phi (\zeta,t)^\top e_d (t) d\zeta, \\
&= \dfrac{1}{2}\int_a^b x_d(t)^\top  \Phi (\zeta,t)\Phi (\zeta,t)^\top e_d (t) d\zeta, \\
&= \dfrac{1}{2}x_d(t)^\top \left( \int_a^b   \Phi (\zeta,t)\Phi (\zeta,t)^\top d\zeta\right) e_d (t) , \\
&= \dfrac{1}{2} x_d(t)^\top E(t) e_d (t) =\dfrac{1}{2} x_d(t)^\top Q(t) x_d (t)  .
\end{split}
\end{equation}
%\begin{equation}\label{DiscretizedHamiltonian}
%\begin{split}
%H_d(t) &= \dfrac{1}{2}\int_a^b x_d(t)^\top \Phi (\zeta,t) \Phi (\zeta,t)^\top x_d (t) d\zeta, \\
%&= \dfrac{1}{2} x_d(t)^\top E(t) x_d (t).
%\end{split}
%\end{equation}

The main properties that we use to develop the time derivative of \eqref{DiscretizedHamiltonian} are:
\begin{align*}
E(t) & = E(t)^\top ,\quad Q(t) = Q(t)^\top,\\
y(t) & = \Phi(a,t)^\top e_d(t), \\
u(t) & = \Phi(b,t)^\top e_d(t),\\
E(t)\dot{x}_d(t) &= \left( \int _a^b \Phi(\zeta,t) \dfrac{\partial \Phi}{\partial \zeta}(\zeta,t)^\top d\zeta \right) e_d(t) \\
& \quad \quad\quad \quad - \left( \int_a^b \Phi(\zeta,t)\dfrac{\partial \Phi}{\partial t}(\zeta,t) ^\top d\zeta\right) x_d(t).
\end{align*}

In the following, we omit time and space dependency and we introduce the following notation: $\Phi_a:= \Phi(a,t)$ and $\Phi_b:= \Phi(b,t)$. Then, we compute the time derivative of the discretized Hamiltonian \eqref{DiscretizedHamiltonian} as follows:

\begin{equation}
\begin{split}
2\dot{H}_d &= {x}_d^\top Q \dot{x}_d + x_d^\top \tfrac{\partial}{\partial t}\left( Qx_d \right)  \\
&= {e}_d^\top E \dot{x}_d + x_d^\top \tfrac{\partial}{\partial t}\left( Qx_d \right)  \\
&= {e}_d^\top \left( F e_d+ B u \right)  + x_d^\top \tfrac{\partial}{\partial t}\left( Qx_d \right),  \\
&=   -{e}_d^\top \Phi_a\Phi_a^\top e_d + {e}_d^\top \Phi_b u -{e}_d^\top \left( \int _a^b \tfrac{\partial \Phi}{\partial \zeta}\Phi^\top d\zeta \right)e_d   \\
& \quad  -{e}_d^\top \left( \int _a^b  \Phi \tfrac{\partial \Phi}{\partial t} ^\top d\zeta \right) x_d    + x_d^\top \tfrac{\partial}{\partial t}\left( Qx_d \right),  \\
&=   -y^2 + u^2 -{e}_d^\top \left( \int _a^b \tfrac{\partial \Phi}{\partial \zeta}\Phi^\top d\zeta \right) e_d   \\
& \quad  -{e}_d^\top \left( \int _a^b  \Phi \tfrac{\partial \Phi}{\partial t} ^\top d\zeta \right)  x_d    + x_d^\top \tfrac{\partial}{\partial t}\left( Qx_d \right).  \\
\end{split}
\end{equation}
Notice that if we are able to show that 
\begin{equation} \label{EqToShow}
\begin{split}
     x_d^\top \tfrac{\partial}{\partial t}\left( Qx_d \right) &= {e}_d^\top \left( \int _a^b \tfrac{\partial \Phi}{\partial \zeta}\Phi^\top d\zeta \right) e_d   \\
     & \quad\quad\quad\quad\quad\quad +{e}_d^\top \left( \int _a^b  \Phi \tfrac{\partial \Phi}{\partial t} ^\top d\zeta \right) x_d 
     \end{split}
\end{equation}
then, 
\[ 2 \dot{H}_d = -y^2 + u^2 \Leftrightarrow \dot{H}_d = \dfrac{1}{2}(u^2-y^2).\]
To do this, we use the following facts:
\begin{equation}
\begin{split}
&{e}_d^\top \left(\int _a^b \tfrac{\partial \Phi}{\partial \zeta}\Phi^\top d\zeta \right) e_d = {e}_d^\top \left( \int _a^b \Phi \tfrac{\partial \Phi}{\partial \zeta}^\top d\zeta \right) e_d, \\
&E\dot{x}_d = \left( \int _a^b \Phi \tfrac{\partial \Phi}{\partial \zeta}^\top d\zeta \right) e_d - \left( \int_a^b \Phi\tfrac{\partial \Phi}{\partial t} ^\top d\zeta \right) x_d.
\end{split}
\end{equation}
Then, \eqref{EqToShow} is equivalently written as:
\begin{equation}
\begin{split}
 x_d^\top & \tfrac{\partial}{\partial t}\left( Qx_d \right) = {e}_d^\top \left( \int _a^b \Phi \tfrac{\partial \Phi}{\partial \zeta}^\top d\zeta \right) e_d      \\
& - {e}_d^\top \left( \int _a^b  \Phi \tfrac{\partial \Phi}{\partial t} ^\top d\zeta\right) x_d  +2 e_d^\top \left(\int _a^b  \Phi \tfrac{\partial \Phi}{\partial t} ^\top d\zeta \right) x_d , \\
&= {e}_d^\top E \dot{x}_d +2 e_d ^\top \left( \int _a^b  \Phi \tfrac{\partial \Phi}{\partial t} ^\top d\zeta \right) x_d , \\
&= {x}_d^\top Q \dot{x}_d +2 e_d ^\top \left( \int _a^b  \Phi \tfrac{\partial \Phi}{\partial t} ^\top d\zeta \right) x_d , \\
 \end{split}
\end{equation}

Finally, if $\mathcal{H}$ is constant for all $\zeta \in \Omega$, one can show that $e_d = \mathcal{H}x_d$.Then, the previous equation becomes:
\begin{equation}
\begin{split}
x_d^\top \tfrac{\partial}{\partial t}\left( E x_d \right) &= x_d^\top E \dot{x}_d + 2x_d^\top \left( \int_a ^b \Phi \tfrac{\partial \Phi}{\partial t}d\zeta \right) x_d, \\
&= x_d^\top E \dot{x}_d + x_d^\top \dfrac{\partial E}{\partial t} x_d, \\
&= x_d^\top \left( E \dot{x}_d + \dfrac{\partial E}{\partial t} x_d \right), \\
\end{split}
\end{equation}
which prove \eqref{EqToShow} and thus $\dot{H}_d = \dfrac{1}{2}(u^2-y^2)$.

\begin{remark}
    Notice that by making $\tfrac{\partial \Phi}{ \partial t} = 0$, the moving mesh becomes the fixed one developped in subsection \ref{subsec:1Dcase} 
\end{remark}

\begin{remark}
If the Hamiltonian density is uniform in space, {\it i.e.,} the transport velocity is the same for all $\zeta \in \Omega$, then the proposed discretization  scheme preserves the structure for any time parametrization of the test functions.
\end{remark}

\section{Numerical Simulations} \label{sec:Simulations}
\begin{figure}[!h]
\begin{center}
\includegraphics[width=0.48\textwidth]{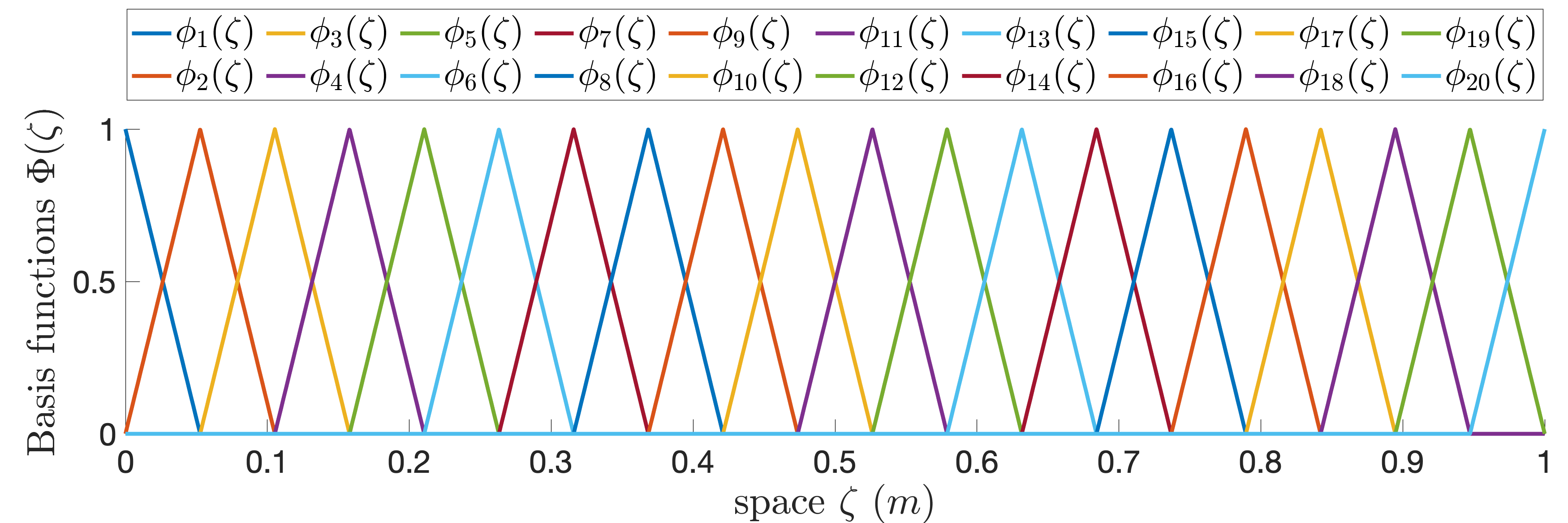}  
\caption{Basis functions $\Phi(\zeta)$ for the fixed mesh.} 
\label{Fig:BasisFunctions_fixed}                    
\end{center}
\end{figure}

\begin{figure}[!h]
\begin{center}
\includegraphics[width=0.48\textwidth]{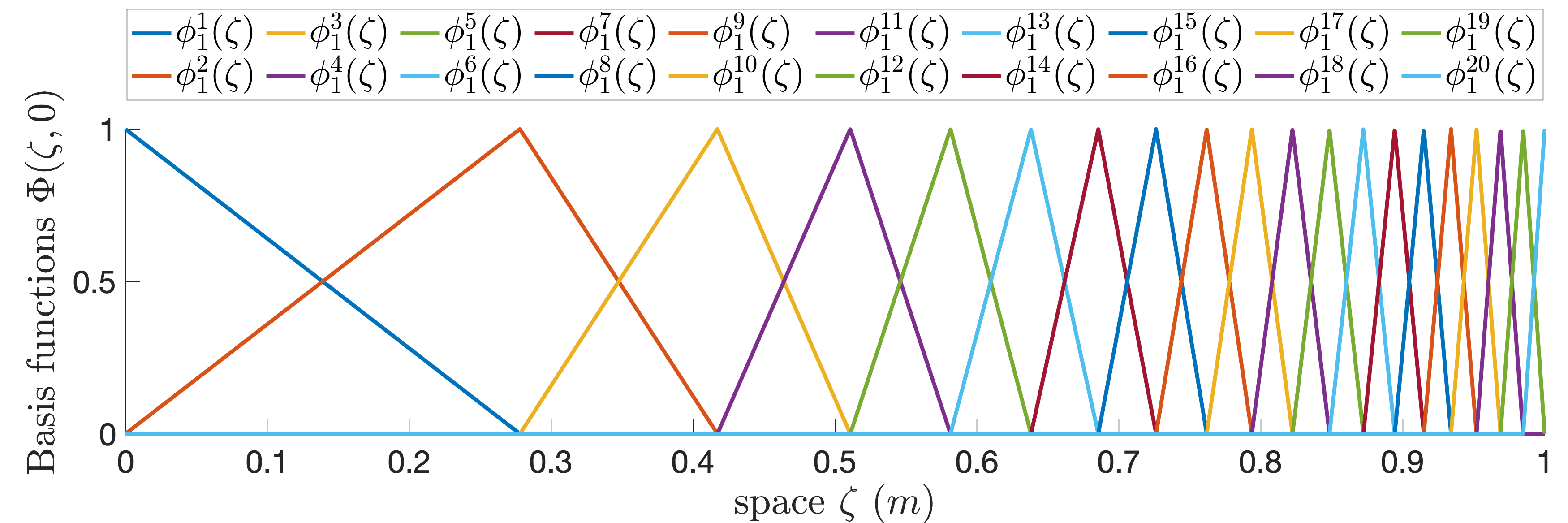}  
\caption{Initial condition of the basis functions $\Phi(\zeta,0)$ for the moving mesh.} 
\label{Fig:BasisFunctionMoving}                      \end{center}
\end{figure}

\begin{figure}[!h]
\begin{center}
\includegraphics[width=0.48\textwidth]{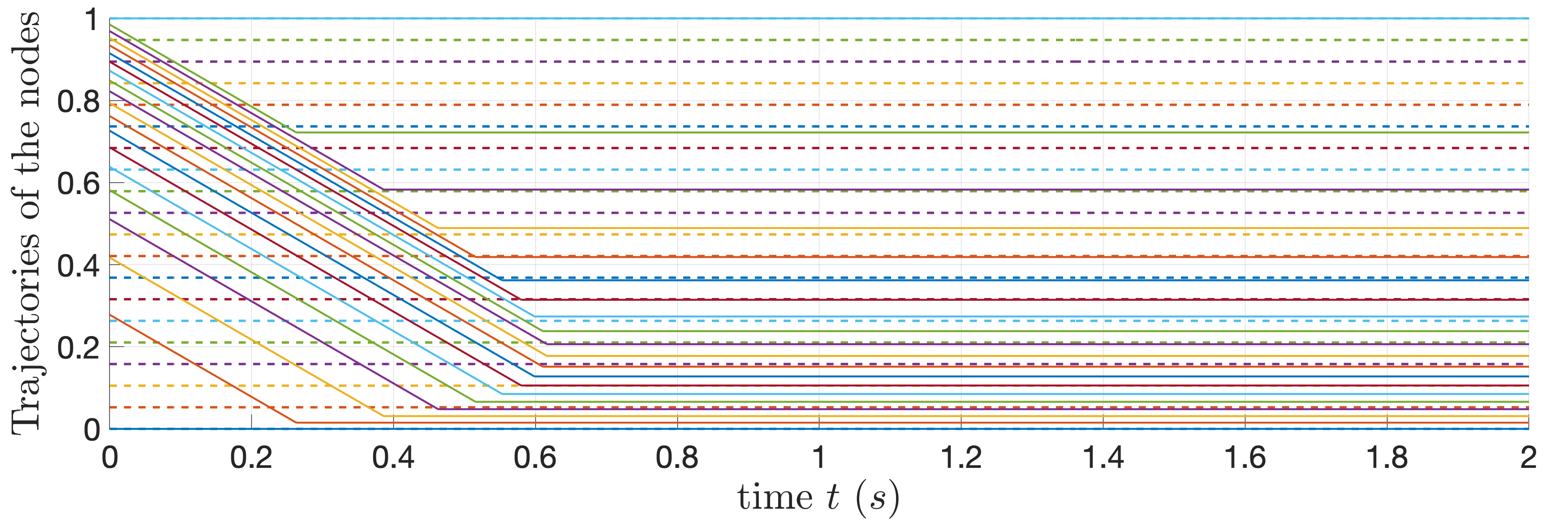}  
\caption{Node trajectories. Dashed lines represent the fixed case and solid ones the moving case.} 
\label{Fig:Trajectories}                      
\end{center}
\end{figure}

We consider the case of the one-dimensional transport equation with constant velocity and we aim at finding a numerical solution using the proposed discretization scheme for a fixed and a moving mesh. The Matlab codes used for this section can be downloaded from a Gitlab repository\footnote{https://gitlab.com/ToledoZucco/moving-mesh-in-1d}. We consider $N = 20$ elements and initial condition $x_0(\zeta) = 2e^{-420(\zeta-0.9)^2}$. In Fig. \ref{Fig:BasisFunctions_fixed}, we show the basis functions used for the case of a fixed mesh in which every node is equidistant with its neighbor nodes. In Fig. \ref{Fig:BasisFunctionMoving}, we show the initial condition of the basis functions for the case of a moving mesh. In this case, the nodes are logarithmically distributed with a bigger concentration at the right side of the spatial domain.
In Fig. \ref{Fig:Trajectories}, we show the nodes trajectories over time. Dashed lines represent the fixed mesh, whereas solid ones represent the moving mesh. For simplicity and since the initial condition has most of its energy concentrated at the right side of the spatial domain, we select a mesh in which the nodes are initially concentrated at the right side. Then, the nodes travel to the left side with a velocity of $1\, (m/s)$ and finally, they end concentrated at the left side. 

\begin{figure}[!h]
\begin{center}
\includegraphics[width=0.48\textwidth]{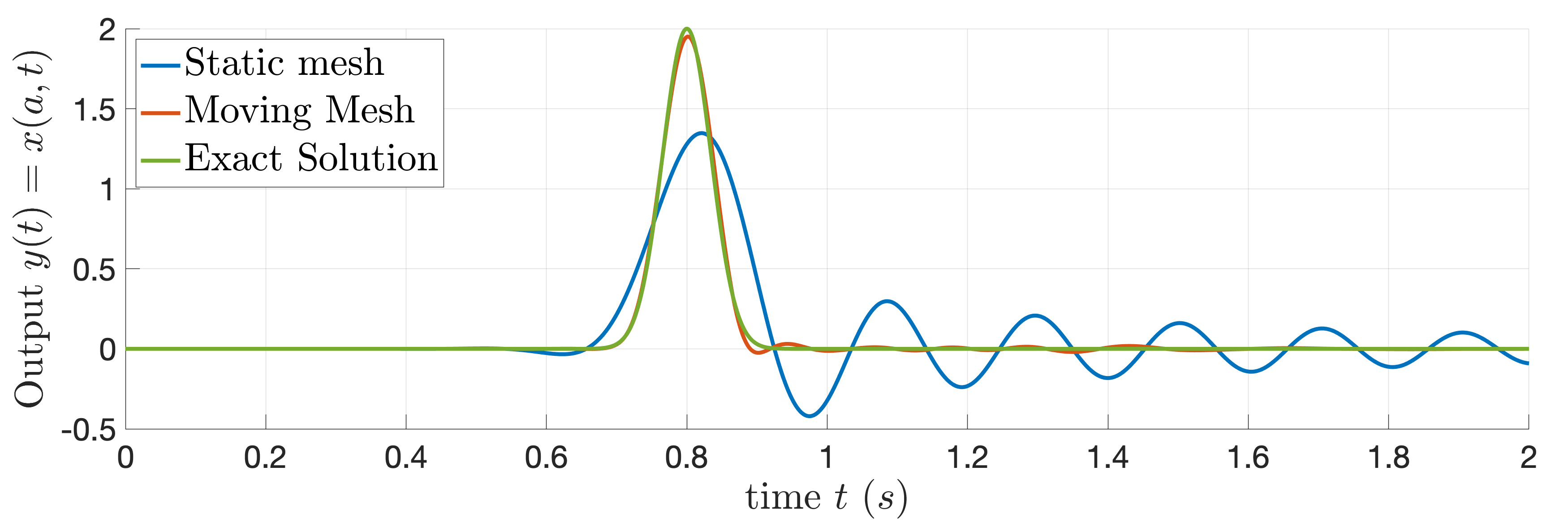}  
\caption{Output of the system $y(t) = x(a,t)$.} 
\label{Fig:Output}                      
\end{center}
\end{figure}

In Fig. \ref{Fig:Output}, we show the exact output of the system and the numerical ones using the fixed and moving mesh scenarios. The solution corresponds to the delayed initial condition. We can see that, in the moving mesh case, the numerical solution is more accurate than in the fixed case. Moreover, less oscillations can be appreciated when using moving meshes. This is a simple example of how moving meshes can reduce oscillations produced when using FEM (typically known as the Gibbs phenomenon).

\begin{figure}[!h]
\begin{center}
\includegraphics[width=0.48\textwidth]{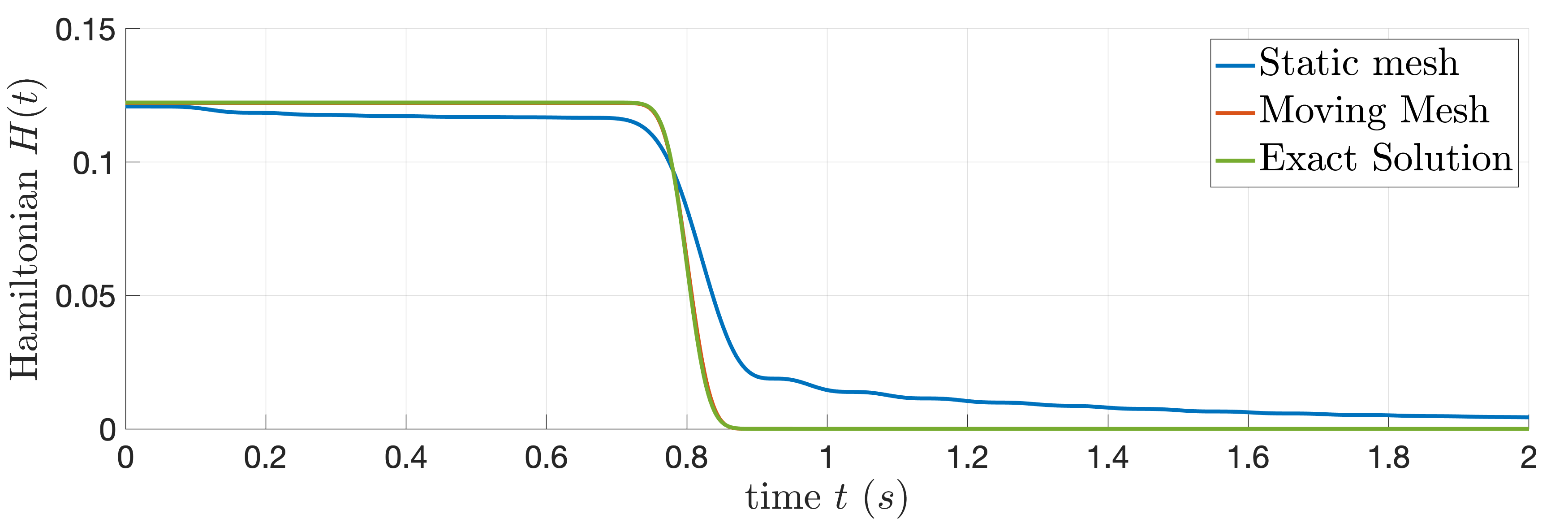}  
\caption{Hamiltonian.} 
\label{Fig:Hamiltonians}                      
\end{center}
\end{figure}

In Fig. \ref{Fig:Hamiltonians}, we show the Hamiltonian over time. We can see that in all cases since the input is zero, there is no increment on the Hamiltonian, and the biggest decrease of the Hamiltonian is around second $t= 0.9\, (s)$, which is the time in which the signal reaches the left boundary.

\begin{figure}[!h]
\begin{center}
\includegraphics[width=0.48\textwidth]{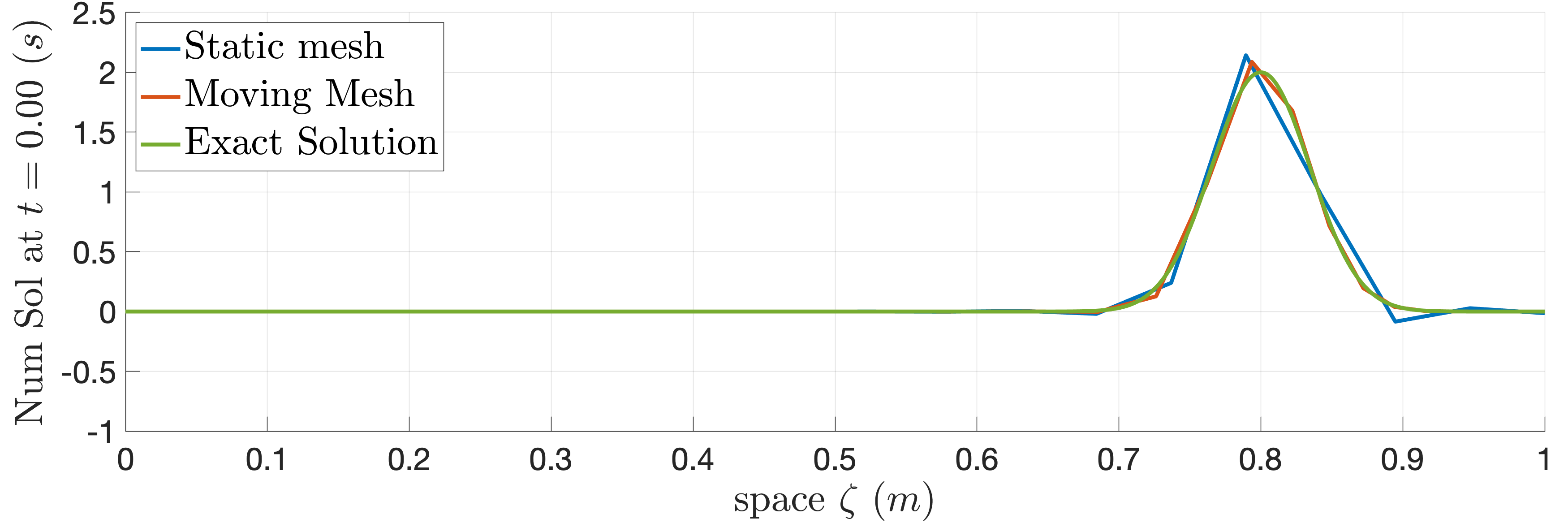}  
\caption{Numerical solution of $x(\zeta,t)$ at $t = 0$ $(s)$.} 
\label{Fig:NumSol1}                      
\end{center}
\end{figure}

\begin{figure}[!h]
\begin{center}
\includegraphics[width=0.48\textwidth]{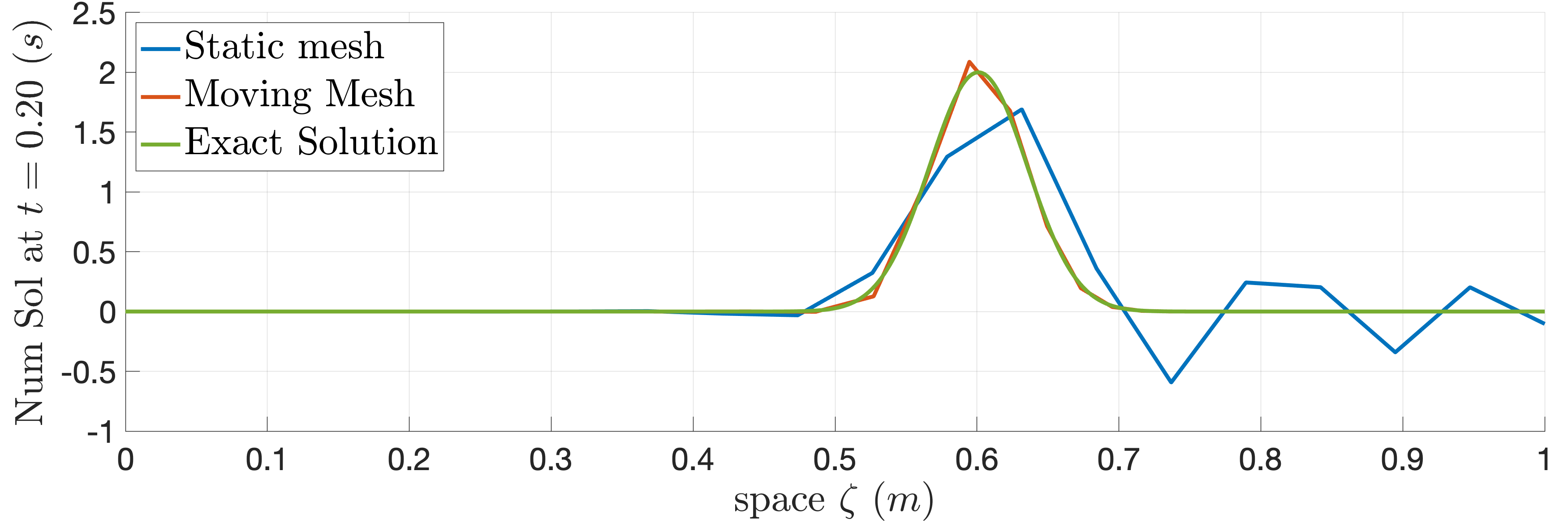}  
\caption{Numerical solution of $x(\zeta,t)$ at $t = 0.2$ $(s)$.} 
\label{Fig:NumSol2}                      
\end{center}
\end{figure}

\begin{figure}[!h]
\begin{center}
\includegraphics[width=0.48\textwidth]{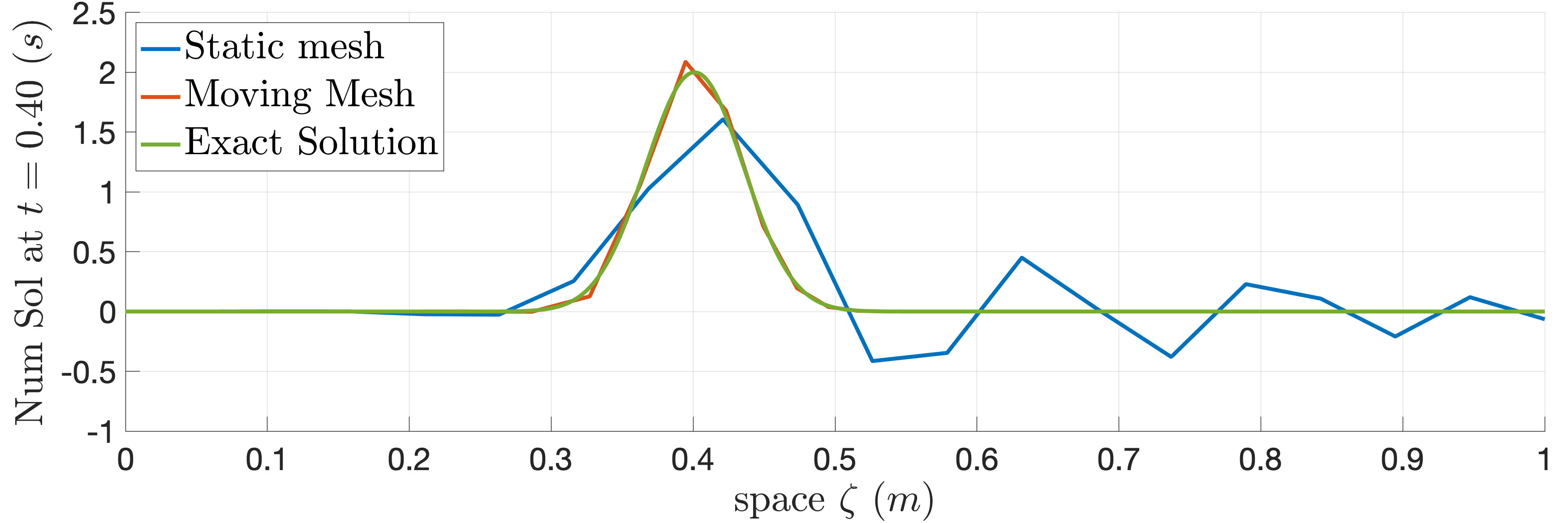}  
\caption{Numerical solution of $x(\zeta,t)$ at $t = 0.4$ $(s)$.} 
\label{Fig:NumSol3}                      
\end{center}
\end{figure}

\begin{figure}[!h]
\begin{center}
\includegraphics[width=0.48\textwidth]{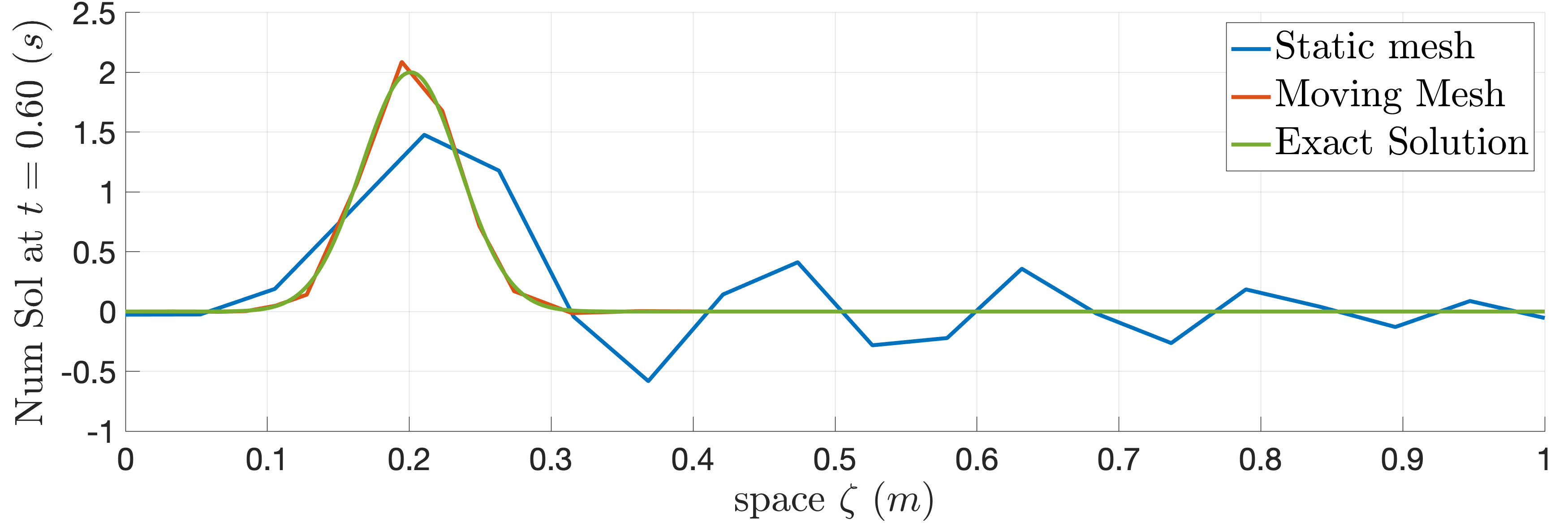}  
\caption{Numerical solution of $x(\zeta,t)$ at $t = 0.6$ $(s)$.} 
\label{Fig:NumSol4}                     
\end{center}
\end{figure}

Finally, in Fig. \ref{Fig:NumSol1}, we show the initial condition $x_0(\zeta)$ and its respective approximation considering a static and moving mesh. The exact solution and the numerical ones are shown at different time steps in Fig. \ref{Fig:NumSol2}, \ref{Fig:NumSol3}, and \ref{Fig:NumSol4}. We can see that the numerical solution is less damaged when using a moving mesh approach.


\section{Conclusions and perspectives}\label{sec:conclusions}
In this paper we have presented a structure-preserving Finite Element Method (FEM) for the transport equation in one- and two-dimensional spatial domains. The main advantages of the proposed method is that the finite-dimensional models preserve the {\it scattering energy preserving} structure of the initial Partial Differential Equation (PDE). Moreover, for the one-dimensional case, we have extended this approach to deal with moving meshes, which can improve the numerical solution behaviour with respect to fixed meshes. This is a first step towards using adaptive meshes whereas preserving the structure. Further investigations on how the nodes should be adapted will be requiered. The extension of the moving mesh approach to {\it passive-energy preserving} is part of a future work in order to deal with waves and beams, for instance. Additionally, the implementation of this discretization scheme in predictive control techniques \citep{Deng2022JournalPredictor} as the well-known Smith predictor is part of the future work of the authors. The fact of preserving the energy structure of the transport phenomena provides a natural Lyapunov function that can be used to guarantee closed-loop stability at the continuos and discrete level \citep{Krstic2008JournalBackstepping}, \citep{BastinCoronBook2016}.  

%\begin{ack}
%\textcolor{red}{Place acknowledgments here.}
%\end{ack}
%\bibliographystyle{apa}

\bibliography{References}             
\end{document}